\pdfoutput=1

\RequirePackage{fix-cm}

\documentclass[leqno, fleqn, 12pt]{article}

\usepackage{latexsym}
\usepackage{amsmath}
\usepackage{amscd}
\usepackage{amssymb}
\usepackage{verbatim}

\usepackage[T1]{fontenc}

\usepackage[upright, widespace]{fourier}

\usepackage{fullpage}

\usepackage{tikz-cd}
\tikzcdset{row sep/mysize6/.initial=4.8em}

\makeatletter
\renewcommand{\@makefntext}[1]{\vspace*{0.5ex}\parindent=0em
\hspace*{-0.4em}
\hbox to 0.4em{\hss\@makefnmark}\hspace*{0.4em}{#1}
}
\makeatother

\newcounter{mysectionnumber}
\newcommand{\mysection}[2]{\setcounter{footnote}{0}
\setcounter{equation}{0}
\setcounter{myparnum}{0}
\refstepcounter{mysectionnumber}
\vspace{24pt}{\Large {\themysectionnumber.} {#1}}\label{#2}\vspace*{24pt}}

\newcounter{myparnum}
\newcounter{mylemmanum}[myparnum]
\newcommand{\myuppar}[1]{\vspace{\medskipamount}\textbf{#1}\hspace*{0.5em}}

\newcommand{\myitpar}[1]{\vspace{\medskipamount}\textit{\textbf{#1}}\hspace*{0.5em}}

\newcounter{myappendnumber}
\newcounter{myaparnum}[myappendnumber]
\newcommand{\proof}{\vspace{\medskipamount}{\textbf{{\emph{Proof}.}}\hspace*{1em}}}

\newcommand{\eproof}{ $\blacksquare$}

\newcommand{\dis}{\displaystyle}

\def\sss{\hspace{0.05em}\ }
\def\dss{\hspace{0.1em}\ }
\def\trs{\hspace{0.15em}\ }
\def\qss{\hspace{0.2em}\ }
\def\pss{\hspace{0.3em}\ }
\def\oss{\hspace{0.4em}\ }

\def\halfff{\hspace*{0.025em}}
\def\fff{\hspace*{0.05em}}
\def\dff{\hspace*{0.1em}}
\def\trf{\hspace*{0.15em}}
\def\qff{\hspace*{0.2em}}
\def\pff{\hspace*{0.3em}}
\def\off{\hspace*{0.4em}}

\newcommand{\nsp}{\hspace*{-0.1em}}
\newcommand{\nnsp}{\hspace*{-0.15em}}

\newcommand{\dnsp}{\hspace*{-0.2em}}

\renewcommand{\leq}{\leqslant}
\renewcommand{\geq}{\geqslant}

\begin{document}

\setlength{\baselineskip}{12pt plus 0pt minus 0pt}
\setlength{\parskip}{12pt plus 0pt minus 0pt}
\setlength{\abovedisplayskip}{12pt plus 0pt minus 0pt}
\setlength{\belowdisplayskip}{12pt plus 0pt minus 0pt}

\newskip\smallskipamount \smallskipamount=3pt plus 0pt minus 0pt
\newskip\medskipamount   \medskipamount  =6pt plus 0pt minus 0pt
\newskip\bigskipamount   \bigskipamount =12pt plus 0pt minus 0pt

\title{Sperner's\qss lemma,\oss
Brouwer's\pss fixed-point\qss theorem,\oss 
and\qss cohomology}
\date{} 
\author{\textnormal{Nikolai\qss V.\qss Ivanov{}\vspace*{6.3ex}}}

\footnotetext{\hspace*{-0.65em}\copyright\oss 
Nikolai\qss V.\qss Ivanov,\oss 2009,\oss 2019.\trs 
Neither the work reported in\dss the present paper\halfff,\qss
nor its preparation were supported by any corporate entity.}

\maketitle

\vspace*{9ex}\vspace*{-6pt}

\renewcommand{\baselinestretch}{1}
\selectfont

\mysection{Introduction}{introduction}

\vspace*{6pt}
The proof\dss of\pss Brouwer's\qss fixed-point\dss theorem based on\dss 
Sperner's lemma\qss \cite{s}\qss is\dss often presented 
as an elementary\sss combinatorial\sss alternative\sss 
to advanced proofs based on algebraic\sss topology\halfff.\pss 
A natural\sss question\dss is\dss to what extent\dss this proof\dss is\dss independent
from\dss the ideas of\dss algebraic\sss topology\halfff,\qss and,\qss
in particular\halfff,\pss to what\sss extent\dss this proof\dss is\dss 
independent\dss from\dss Brouwer's\dss own proof\halfff.\oss
The\sss latter\dss is\dss based on\dss the degree\qss theory\halfff,\qss 
i.e.\qss on a fragment\sss of\dss algebraic\sss topology\sss developed\dss by\dss Brouwer\halfff.\pss 
After\sss the author discovered\qss \cite{i}\qss that\dss the 
famous analytic proof\dss of\pss Brouwer's\qss theorem 
due to\dss Dunford\dss and\dss Schwartz\dss is\dss nothing else but\dss 
the usual\dss topological\dss proof\dss in disguise,\qss 
it\dss was only\dss natural\dss to suspect\dss that\dss the same\dss is\dss 
true for the proof\dss based on\qss Sperner's\trs lemma.\pss 
This suspicion\dss turned out\dss to be correct\halfff,\qss 
and\dss the goal of\dss this note\dss is\dss to uncover\dss 
the standard\dss topology\dss hidden\dss in\dss this proof\halfff.\qss

In fact,\qss the\sss two situations are very similar.\qss 
Dunford--Schwartz\dss proof\dss can be considered as 
a cochain-level\dss version of\dss the standard proof\dss 
based on\dss de Rham\dss cohomology\dss theory\qss 
(in\dss the\dss de Rham\dss theory cochains are nothing\sss else but\dss differential forms),\qss 
written\dss in the\sss language of\dss elementary\dss multivariable calculus.\qss 
Similarly\halfff,\qss the combinatorial proof\dss of\pss Sperner's\trs 
lemma can be considered as a cochain-level version 
of\dss a standard argument\dss based on simplicial\sss cohomology\dss theory\halfff,\qss 
written in a combinatorial language.\qss
Of\dss course,\qss in\dss this version\sss one needs to use simplicial\sss cochains 
instead of\qss de Rahm\dss cochains\qss (i.e.,\qss of\dss differential forms).\qss
 
Both\sss these alternatives\sss to algebraic\sss topology\dss 
turn out\dss to be not\dss homological,\qss 
but\sss cohomological\dss proofs in disguise.\qss 
This\dss is\dss somewhat\sss surprising\fff;\qss 
usually\dss homology theory\dss is\dss considered\dss 
to be more intuitive than cohomology theory\halfff.\qss
The most\dss likely reason\dss is\dss that\sss cohomology\dss theory\dss is\dss a\qss
\emph{contravariant}\pss functor\fff:\oss
preimages under maps behave better\dss than images.

Contrary\dss to a widespread\dss belief\halfff,\oss
Sperner\dss did\dss not\dss proved\dss his lemma in order\sss to
provide an elementary\dss proof\dss of\pss Brouwer's\trs fixed\dss point\dss theorem.\oss
He was interested in other\sss theorems of\qss Brouwer\qss
(the invariance of\dss dimension and\sss the invariance of\dss domains),\oss
and\qss Brouwer's\dss fixed\dss point\dss theorem\dss is\dss not\sss even\dss
mentioned in his paper\qss \cite{s}.\oss
The standard deduction of\qss Brouwer's\dss theorem 
from\dss Sperner's\dss lemma\dss is\dss 
based on an argument of\qss
Knaster--Kuratowski--Mazurkiewich\qss \cite{kkm}.\oss
At\dss the first\sss sight\dss this argument\dss is\dss quite different\dss
from arguments used\sss in proofs
of\qss Brouwer's\dss fixed\dss point\dss theorem\dss 
based on algebraic\sss topology\halfff.\oss 

These proofs of\pss Brouwer's\qss theorem are usually\dss based on a 
construction of\dss a retraction of\dss a disc onto its boundary\sss starting\dss
from a fixed point free self-map of\dss a disc.\qss
There are no such retractions by the\qss \emph{no-retraction\dss theorem}\qss and\qss
Brouwer's\dss theorem follows.\qss
It\dss turns out\dss that\dss this construction of\dss a retraction from a fixed point free map
underlies also\sss the\qss Knaster--Kuratowski--Mazurkiewich\qss argument\halfff.\oss
The no-retraction\dss theorem itself can be proved\dss by a modification of\dss this 
argument\sss suggested\dss by\dss the notion of\qss \emph{simplicial approximations}.\qss

\myuppar{Outline of\dss the paper.} 
The rest of\dss the paper\dss is\dss divided into\sss two parts.\pss 
Section\qss \ref{sperner-cohomology}\qss is\dss devoted to\qss Sperner's\trs lemma.\pss 
We start with\dss Sperner's\dss own proof\dss phrased\dss in a\dss geometric\sss language 
and\dss then present\sss a cohomological\dss proof\halfff.\pss 
In order to compare these two proofs we rewrite\sss 
the cohomological\dss proof\dss in\dss terms of\dss cochains.\pss 
The readers not\dss familiar with the cohomology theory may skip the cohomological\dss proof\dss
and\dss proceed directly\dss to this cochain-level\dss proof\halfff.\pss
The\sss latter may\dss be thought\sss of\dss as a\sss 
geometric realization of\qss Sperner's\dss arguments.\qss

The version of\dss cohomology theory\dss most suitable for discussing\qss
Sperner's\dss proof\dss is\dss the simplicial\sss cohomology\dss theory\halfff.\pss 
Some familiarity\dss with\dss the cohomology\dss theory on\dss 
the part\sss of\dss the reader will\dss help,\qss 
but\dss the author hopes that\dss the 
cochain-level\dss proof\dss
will\dss be accessible\sss to all\dss readers comfortable with\dss Sperner's\dss lemma\dss
and\dss linear algebra methods in combinatorics.\pss 

Section\qss \ref{deduction}\qss is devoted to\dss Brouwer's\dss fixed-point\dss theorem.\pss 
We present\qss Knaster--Kuratowski--Mazurkiewich\qss deduction of\qss 
Brouwer's\trs theorem\dss from\qss Sperner's\trs lemma 
and\dss then explain\dss how\dss the no-retraction\dss theorem 
and its proof\dss are,\qss 
in fact,\qss hidden\dss in\dss this deduction.\oss
In\dss particular\halfff,\pss 
one can use\qss Sperner's\trs lemma\dss in order\dss to prove\dss 
the\sss no-retraction\dss theorem
and\dss then complete\sss the proof\dss of\pss Brouwer's\trs theorem\dss 
by\dss the well\dss known elementary\dss 
geometric argument\halfff.\oss

\myuppar{Acknowledgments.} 
I\dss am\dss grateful\dss to\qss F.\qss Petrov\qss and\dss the\sss late\qss 
A.\qss Zelevinsky\qss for\dss their interest\dss 
in cohomological\dss interpretation of\pss Sperner's\trs lemma,\oss 
which stimulated\dss me to write\sss this paper\halfff.\oss 
I\dss am\dss also grateful\dss to\qss
M.\qss Prokhorova\qss for careful\dss reading of\dss the first\dss 
version of\dss this paper\halfff,\pss
for\dss pointing\sss out\dss that\dss the standard construction of\dss 
the above-mentioned\dss retraction\dss may\dss be not\dss well-defined\dss 
for\sss a simplex\halfff,\oss 
and\dss for continuing\dss interest\dss in\sss my\dss work\halfff.\pss

\newpage
\mysection{Sperner's\qss lemma\qss and\qss cohomology}{sperner-cohomology}

\vspace*{5.25pt}
Let\sss $\Delta$\sss be an $n${\dnsp}-dimensional simplex with the vertices\qss 
$v_0\dff,\dff v_1,\dff\ldots\dff,\dff v_{n}$\nsp.\pss 
Let\sss $\Delta_{\dff i}$\sss be\sss the face of\dss $\Delta$\dss opposite\sss to\sss the vertex\dss $v_i$\nsp.\oss 
Its\dss vertices are all\dss the vertices\qss 
$v_0\dff,\dff v_1,\dff\ldots\dff, v_{n}$\qss except $v_i$\nsp.\oss 
The boundary\sss $\partial\dff \Delta$\sss is\dss equal\dss to\sss 
the union of\trs the\dss $n\qff +\qff 1$\dss faces\dss $\Delta_{\dff i}$\nsp.\qss 
Suppose that the simplex\dss $\Delta$\dss 
is subdivided into smaller simplices forming a simplicial complex\dss $S$\nnsp.\qss\vspace{-1pt}

\myuppar{Sperner's\dss lemma.} 
\emph{Suppose that vertices of\trs $S$\sss 
are labeled by the numbers\qss $1\dff,\dff 2\dff,\dff\ldots\dff,\dff n+1$\qss 
in such a way that if\dss a vertex\sss $v$\sss belongs to a face\sss $\Delta_{\dff i}$\nnsp,\qss 
then the label of\trs $v$\sss is not equal to\sss $i$\nnsp.\oss 
Then the number of\trs $n${\dnsp}-dimensional simplices of\qss $S$\dss 
such that the set of\qss labels of\qss their vertices\dss is\dss equal\dss to\qss  
$\{\qff 1\dff,\dff 2\dff,\dff\ldots\dff,\dff n\qff +\qff 1 \qff\}$\qss is\dss odd.\oss
In particular,\pss there\dss is\dss at\dss least\sss one such simplex.\qss}\vspace{-1pt}

\myuppar{Geometric interpretation.}
The\sss labelings\sss in\qss Sperner's\dss lemma admit\sss a geometric interpretation.\pss 
Namely,\qss if\dss a vertex\dss $v$\dss is\dss labeled\dss by $i$\nnsp,\oss 
let\sss us set\qss $\varphi\dff(\fff v\trf)\qff =\qff v_i$\nnsp.\oss 
This defines a map $\varphi$ from\sss the set\sss of vertices of\sss $S$ 
to\sss the set\dss of vertices of\sss $\Delta$\nnsp.\qss 
Since $\Delta$ is\dss a\sss simplex and\dss therefore every set of\dss vertices of\sss $\Delta$ 
is\dss a set\sss of\dss vertices of\dss a subsimplex,\qss 
$\varphi$ defines a\qss \emph{simplicial map}\qss 
$S\qff \longrightarrow\qff \Delta$\nnsp,\qss which we will\sss also denote by $\varphi$\nnsp.\pss 
Clearly,\qss an $n${\dnsp}-dimensional\sss simplex\dss has\qss
$\{\qff 1\dff,\dff 2\dff,\dff\ldots\dff,\dff n\qff +\qff 1 \qff\}$\qss 
as\sss the set of\trs labels of\dss its vertices\dss  
if\trs and\dss only\trs if\dss $\varphi$ maps it\sss onto $\Delta$\nnsp.\qss 
In\dss this language\sss the assumption of\qss Sperner's\trs 
lemma means\sss that $\varphi$ maps 
simplices of\sss $S$ contained\dss in\sss a face\dss $\Delta_{\dff i}$\dss 
of\dss $\Delta$\dss into\dss $\Delta_{\dff i}$\nnsp,\oss 
and\dss the first\dss part\sss of\dss the conclusion means\sss that $\varphi$ 
maps an odd\dss number of\dss simplices onto $\Delta$\nnsp.\qss
The second part\dss is\dss deduced\dss from\dss the first\sss one by\sss applying one of\dss
the most\dss fundamental\dss principles of\dss mathematics:\qss 
\emph{an odd\dss number\dss is\dss not\dss equal\dss to zero.}

In\dss their classical\dss treatise\qss \cite{ah}\qss Alexandroff\dss and\dss Hopf\qss state\qss 
Sperner's\trs lemma in this language 
and\dss use algebraic topology\dss to prove it\dss
without\sss even\dss mentioning\dss the combinatorial approach.\oss
For quite a while the author\sss naively\dss thought\dss that\sss he\sss rediscovered
this geometric interpretation.\qss
But\dss nowadays\sss the idea of\dss turning sets into simplicial complexes is well-established,\oss
and one can\dss trace its roots to the notion of\dss the nerve of\dss a system of\dss
sets,\qss introduced\dss by\dss Alexandroff\qss \cite{a-nerve}.\oss
The\qss ``rediscovery''\qss was,\oss albeit\dss indirectly\halfff,\oss 
based on\dss the original\dss discovery.\oss\vspace{-1pt}

\myitpar{Combinatorial proof\dss of\qss Sperner's\qss lemma.} 
Now we will\sss present a\sss translation of\qss Sperner's\dss proof\sss
in the geometric language.\oss
Along the way we introduce some notations used later.\qss 
The numbers\qss $e\fff,\dff f\fff,\dff g\fff,\dff h$\qss 
introduced in\dss this proof\dss have exactly\dss the same meaning as 
in\dss Sperner's\dss paper\qss \cite{s}\qss and\dss 
play a crucial\dss role also in\dss the cochains-based\dss proof\trs later\halfff.\qss\vspace{-1pt}

Let us use induction by $n$\nnsp.\pss 
The result is trivially true for\qss $n\qff =\qff 0$\nnsp.\pss
Suppose that\qss $n\qff >\qff 0$\nnsp.\pss
Let us consider one of\dss the faces of\dss $\Delta$\nnsp,\qss
say,\qss the face\dss $\Delta_{\dff n\dff +\dff 1}$\nnsp.\pss
It is an $(\fff n\dff -\dff 1\fff)$\dnsp-dimensional simplex.\qss
The subdivision\dss $S$\dss defines a subdivision\dss $S_{\dff n\dff +\dff 1}$\dss 
of\dss $\Delta_{\dff n\dff +\dff 1}$\dss
and the restriction of\dss $\varphi$\dss to\dss $S_{\dff n\dff +\dff 1}$\dss is a simplicial map\qss 
$S_{\dff n\dff +\dff 1}\qff \longrightarrow\qff \Delta$\nnsp.\oss
By the assumptions of\qss Sperner's\dss lemma\sss the image of\dss this map
is contained in $\Delta_{\dff n\dff +\dff 1}$  
and hence $\varphi$ 
defines a simplicial map\qss 
$S_{\dff n\dff +\dff 1}\qff \longrightarrow\qff \Delta_{\dff n\dff +\dff 1}$\nnsp.\oss
The assumptions of\qss Sperner's\dss lemma obviously hold\dss for this map,\qss 
and\dss by\dss the inductive assumption $\varphi$ 
maps an odd number of\dss simplices of\qss 
$S_{\dff n\dff +\dff 1}$ onto\dss $\Delta_{\dff n\dff +\dff 1}$\nnsp.\qss 
It remains to deduce from this that $\varphi$ maps an odd number of simplices of\dss 
$S$\dss onto\dss $\Delta$\nnsp.\oss 
This is\dss the main part of\dss the proof.\qss

If\sss $\sigma$ is a simplex of\dss $S$\sss contained in the boundary $\partial\dff \Delta$
and
$\varphi\dff(\dff \sigma\dff)
\off =\off 
\Delta_{\dff n\dff +\dff 1}$\nnsp,\oss
then $\sigma$ is contained\dss in $\Delta_{\dff n\dff +\dff 1}$\nnsp.\oss
Indeed,\oss otherwise $\varphi\dff(\dff \sigma\dff)$ 
would be contained in some other face by the assumptions
of\qss Sperner's\dss Lemma.\oss
Let\qss $\sigma_{\dff 1}\dff,\pff \ldots\dff,\pff \sigma_{\dff h}$\qss
be the\sss full\dss list of\dss such simplices.\oss
All of\dss them are contained in\dss $\Delta_{\dff n\dff +\dff 1}$\nnsp.\oss
Let\qss $\tau_{\dff 1}\dff,\pff\ldots\dff,\pff \tau_{\dff g}$\qss 
be the\sss full\dss list of other $(\fff n\dff -\dff 1)$\dnsp-dimensional 
simplices of\dss $S$\sss mapped by $\varphi$ onto $\Delta_{\dff n\dff +\dff 1}$\nnsp.\qss 
They are contained in the interior of\dss $\Delta$\nnsp.

Let us consider $n${\dnsp}-dimensional simplices of\dss $S$\nnsp.\oss 
We are interested only in $n$\dnsp-dimensional simplices $\sigma$ such\dss that\dss the image 
$\varphi\dff(\dff \sigma\dff)$ is equal either\dss to
$\Delta$ or\dss to $\Delta_{\dff n\dff +\dff 1}$\nnsp.\qss
Clearly\halfff,\oss on such simplices have a face mapped\dss by $\varphi$
onto $\Delta_{\dff n\dff +\dff 1}$\nnsp.\qss
Let\qss $\rho_{\dff 1}\dff,\pff \ldots\dff,\pff \rho_{\dff e}$\qss be the full\sss list of 
$n$\dnsp-dimensional simplices  of\dss $S$ having $\Delta$ as\sss the image
under $\varphi$\nnsp.\oss
Each of\dss these simplices has exactly one face mapped\dss by $\varphi$
onto $\Delta_{\dff n\dff +\dff 1}$\nnsp.\qss

Let $f$ be the number of\dss $n${\dnsp}-dimensional simplices $\sigma$ such\dss that\dss
$\varphi\dff(\dff \sigma\dff)\off =\off \Delta_{\dff n\dff +\dff 1}$\nnsp.\oss
We claim that every such simplex has exactly two faces
mapped\sss by $\varphi$ onto $\Delta_{\dff n\dff +\dff 1}$\nnsp.\qss
Indeed,\oss such simplex $\sigma$ obviously\sss has
a proper face $\tau$ mapped by $\varphi$ onto $\Delta_{\dff n\dff +\dff 1}$\nnsp.\qss
If $v$ is the vertex of $\sigma$ not belonging to this face,\pss 
then\qss $\varphi\dff(\fff v\dff)\qff =\qff v_i$\qss for some\qss 
$i\qff \leqslant\qff n$\nnsp.\qss 
Let $w$ be the vertex of\dss $\tau$ mapped\dss by $\varphi$ to $v_i$\nnsp.\qss 
Clearly,\pss if\trs we replace $w$ by $v$ in $\tau$\nnsp,\pss 
we will\dss get another face $\tau'$ of\dss $\sigma$ mapped\dss 
by $\varphi$ onto\dss $\Delta_{\dff n\dff +\dff 1}$\nnsp,\qss
and no other face of\dss $\sigma$ is\dss mapped onto\dss $\Delta_{\dff n\dff +\dff 1}$\nnsp.\oss 
This proves our claim.\oss 

\emph{Double counting.}\oss
Now we are ready for\dss the most acclaimed part of\dss the proof,\pss
a double counting argument.\pss
Let us count in two ways the number of pairs $(\sigma,\dff \tau\dff)$ 
such that $\tau$\dss is an 
$(\fff n\dff -\dff 1\fff)${\dnsp}-dimensional face of an 
$n$\dnsp-dimensional simplex $\sigma$ and $\varphi$ maps 
$\tau$ onto $\Delta_{\dff n\dff +\dff 1}$\nnsp.\pss 
If\dss we count\dss the simplices $\tau$ first,\oss 
we get\qss $2\dff g\dff +\dff h$\qss 
as the number of\dss such pairs,\qss 
because every $(\fff n\dff -\dff 1\fff)$\dnsp-dimensional simplex in the interior of\dss $\Delta$ 
is a face of\dss exactly two $n$-dimensional simplices,\qss 
and every $(\fff n\dff -\dff 1\fff)$\dnsp-dimensional simplex in the boundary $\partial\dff \Delta$ 
is a face of exactly one  $n$\dnsp-dimensional simplex.\qss 
If\dss we count\dss the simplices $\sigma$ first,\qss 
then only the simplices having as\sss the image either $\Delta$ or $\Delta_{\dff n\dff +\dff 1}$ matter and
we get\qss $e\dff +\dff 2\fff f$\qss
as the numbers of pairs.\qss 
Therefore\vspace{3pt}
\begin{equation}
\label{sperner}
\quad
2\dff g\qff +\qff h\off =\off e\qff +\qff 2\fff f\qff,
\end{equation}

\vspace{-9pt}
and hence $e$ is odd\qss if\sss $h$\sss is.\pss 
But $h$ is odd\dss by the inductive assumption 
and hence $e$ is also odd.\pss 
Since $e$ is the number of  
simplices mapped onto $\Delta$\nnsp,\qss 
this completes the proof. \eproof

\myitpar{Cohomological proof of\dss Sperner's\dss lemma.} 
The readers uncomfortable with cohomology\dss theory may skip this proof\halfff.\oss
We will use\sss the simplicial cohomology\sss theory with coefficients in the 
field $\mathbb{F}_{\dff 2}$\nnsp.\qss
Since only in the distinction between even and odd numbers matters,\qss
this is the most natural choice of coefficients.\oss
We will consider the boundary $\partial\dff \Delta$ as a simplicial complex 
having\qss $\Delta_1\fff,\dff \ldots\fff,\dff\Delta_{\dff n\dff +\dff 1}$ 
as the top-dimensional simplices. 
Let 
$\partial\dff S$ be the simplicial complex consisting of simplices of 
$S$ contained in\dss the\qss ({\fff}geometric)\qss boundary of  $\Delta$\nnsp.\qss 
Clearly,\qss $\varphi\colon S\qff \longrightarrow\qff \Delta$\qss 
induces a simplicial map\qss 
$\partial\dff S\qff \longrightarrow\qff \partial\dff \Delta$\nnsp,\qss 
which we will denote by\dss $\varphi_{\dff\partial}$\nnsp.\qss
As in the combinatorial proof,\qss 
we will use an induction by $n$\nnsp,\qss 
the case\qss $n\qff =\qff 0$\qss being trivial.\qss 
Suppose\sss that\qss $n\qff >\qff 0$\qss 
and\dss consider\dss the following\sss commutative diagram.\vspace*{6pt}
\begin{equation}
\label{hdiag}
\quad
\begin{tikzcd}[column sep=huge, row sep=mysize6]\dis
H^{\dff n\dff -\dff 1}\dff (\trf \partial\dff \Delta \trf) 
\arrow[d, "\dis \qff\varphi_{\dff \partial}^*"]
\arrow[r, "\dis \delta"]
& 
H^{\dff n}\dff (\trf \Delta,\qff \partial\dff \Delta\trf)
\arrow[d, "\dis \trf\varphi^*"] 
\\ 
H^{\dff n\dff -\dff 1}\dff (\trf \partial\dff S\trf)
\arrow[r, "\dis \delta"]
& 
H^{\dff n}\dff (\trf S\fff,\qff \partial\dff S\trf)
\off.
\end{tikzcd}
\end{equation}

\vspace{-6pt}
Here the maps $\delta$ are the connecting homomorphisms 
in the cohomological sequences of\dss pairs\qss 
$(\dff S,\dff \partial\dff S \dff)$\qss and\qss
$(\dff \Delta\fff,\dff \partial\dff \Delta \dff)$\nnsp.\qss 
Since the cohomology groups of\dss both $S$ and $\Delta$ are trivial,\qss 
the exactness of the cohomological sequences of\dss pairs 
$(\dff S,\dff \partial\dff S \dff)$ and 
$(\dff \Delta\fff,\dff \partial\dff \Delta \dff)$ implies that both maps 
$\delta$ are isomorphisms.\qss 
In fact,\qss
all cohomology groups in\qss (\ref{hdiag})\qss are isomorphic to $\mathbb{F}_{\dff 2}$\nnsp.\qss
 
The map $\varphi_{\dff \partial}^*$ is the multiplication by the 
${\rm mod}\qff 2$ degree of $\varphi_{\dff \partial}$\nnsp.\qss 
This degree is equal to the number of 
$(\fff n\dff -\dff 1\fff)$\dnsp-simplices mapped by 
$\varphi_{\dff \partial}$ onto any $(\fff n\dff -\dff 1\fff)$\dnsp-simplex of\dss 
$\partial\Delta$\nnsp,\qss 
for example,\qss onto $\Delta_{\dff n\dff +\dff 1}$\nsp.\qss 
In\sss the notations of\dss the combinatorial proof,\qss 
there are $h$ such simplices and\dss hence\sss
the ${\rm mod}\pff 2$ degree 
of\dss $\varphi_{\dff \partial}$ is equal\dss to $h\off {\rm mod}\qff 2$\nnsp.\oss
Therefore\dss 
$\varphi_{\dff \partial}^*$\dss is\sss equal\dss to the multiplication\dss by $h$\nnsp.\qss 
Similarly,\qss 
$\varphi^*$\dss is\dss equal\dss to the multiplication\dss by $e$\nnsp,\qss 
the number of simplices of\dss $S$ mapped onto $\Delta$\nnsp.\qss
Since the maps $\delta$ are isomorphisms,\qss 
the commutativity of\pss (\ref{hdiag})\qss implies that\dss 
$e\qff \equiv\qff h\off {\rm mod}\pff 2$\nnsp.\qss 
Since $h$ is\dss assumed to be odd,\qss 
$e$ is\dss odd\sss also.\oss \eproof

\myuppar{Simplicial cochains.} 
Our next\dss goal\dss is\dss to present a cochain-level version of\dss 
the above cohomological proof\halfff.\qss 
Let\dss us\sss begin with introducing\dss the notions involved.\oss
Let $X$ be a finite simplicial complex and\dss let $A$ be a subcomplex of\dss $X$\nnsp.\qss 
An \emph{$n$\dnsp-dimensional cochain}\qss of\dss $X$ is defined as 
a\sss function on the set of\dss simplices of\dss $X$ of\dss dimension $n$ 
with values in\sss the field\dss $\mathbb{F}_{\dff 2}$\nsp.\qss
Such cochains form a vector space over $\mathbb{F}_{\dff 2}$\nsp,\qss 
denoted\dss by $C^{\dff n}\dff(\trf X\trf)$\nnsp.\qss 
A cochain belongs to the subspace of\qss \emph{relative cochains}\dss 
$C^{\dff n}\dff(\trf X\fff,\qff A\trf)$\dss 
if\qss it\dss is\dss equal\dss to $0$ on all\dss simplices of\dss $A$\nnsp.\qss

To every $n$\dnsp-simplex $\tau$ corresponds a cochain taking the value 
$1$ on $\tau$ and\dss the value $0$ on all other simplices.\qss 
By an abuse of notations,\qss 
we will denote this cochain also by $\tau$\nnsp.\qss 
Such cochains form a basis of $C^{\dff n}\dff(\trf X\trf)$, 
and the cochains corresponding to $n${\dnsp}-simplices not contained in $A$ 
form a basis of $C^{\dff n}\dff(\trf X\fff,\qff A\trf)$\nnsp.\qss
Our\sss abuse of\dss notations allows us to consider cochains as 
formal sums of\dss simplices.\qss 
From this point of\dss view,\qss 
a cochain belongs to\dss $C^{\dff n}\dff(\trf X\fff,\qff A\trf)$\dss 
if\trs and\dss only\trs if\trs
it\dss does not involve simplices of\dss $A$\nnsp.\oss 
With\dss these conventions,\qss the\qss \emph{coboundary map}\vspace{3pt} 
\[
\quad
\partial^* \colon\qff C^{\dff n}\dff(\trf X\trf)
\off \longrightarrow\off 
C^{\dff n\dff +\dff 1}\dff(\trf X\trf)
\]

\vspace{-9pt} 
can be defined as follows.\qss 
If\dss $\tau$ is\dss an $n${\dnsp}-simplex considered as a cochain,\qss 
then $\partial^*(\dff\tau\dff)$ 
is\dss the sum of\dss all $(\fff n\dff +\dff 1\fff)$\dnsp-dimensional 
simplices of\dss $X$ having $\tau$ as a face.\qss
This\sss map $\partial^*$ is\dss extended 
to the whole space $C^{\dff n}\dff(\trf X\trf)$ by\dss linearity.\oss 
Clearly\halfff,\qss
$\partial^*$ maps\dss 
$C^{\dff n}\dff(\trf X\fff,\qff A\trf)$\dss to\dss 
$C^{\dff n\dff +\dff 1}\dff(\trf X\fff,\qff A\trf)$\nnsp.\qss

\myuppar{Simplicial cohomology.}
A cochain\qss $\alpha\qff \in\qff C^{\dff n}\dff(\trf X\trf)$\qss 
is called a\dss \emph{cocycle}\pss if\pss 
$\partial^*\dff(\dff \alpha\dff)\off =\off 0$\qss
and\dss is called a\dss \emph{coboundary}\qss if\trs it\dss
belongs\dss to\sss the image of\qss
$\partial^* \colon\dff C^{\dff n\dff -\dff 1}\dff(\trf X\trf)
\qff \longrightarrow\qff 
C^{\dff n}\dff(\trf X\trf)$\nnsp.\qss
The\dss \emph{relative cocycles}\qss and\dss \emph{coboundaries}\qss
are defined similarly\halfff.\qss
The\dss \emph{cohomology\dss group}\dss
$H^{\dff n}\dff(\trf X\trf)$\dss
is\dss defined as the quotient of\dss the\dss $\mathbb{F}_{\dff 2}$\dnsp-vector space
of\dss $n$\dnsp-dimensional cocycles by\dss the subspace of\dss 
$n${\dnsp}-dimensional coboundaries.\qss
The\dss \emph{relative cohomology group}\qss
$H^{\dff n}\dff(\trf X\fff,\qff A\trf)$\dss
is defined similarly.\qss
So,\qss a cocycle\qss 
$\alpha\qff \in\qff C^{\dff n\dff -\dff 1}\dff(\trf X\trf)$\qss 
defines an element of\dss 
$H^{\dff n\dff -\dff 1}\dff(\trf X\trf)$\nnsp,\qss 
the\dss \emph{cohomology class}\dss of\dss $\alpha$\nnsp,\qss 
and a cocycle\qss
$\alpha\qff \in\qff C^{\dff n\dff -\dff 1}\dff(\trf X\fff,\qff A \trf)$\dss
an element of\dss $H^{\dff n\dff -\dff 1}\dff(\trf X\fff,\qff A\trf)$\nnsp.\oss
The\dss \emph{connecting map}\vspace{3.8pt}
\[
\quad
\delta\qff \colon\qff H^{\dff n\dff -\dff 1}\dff(\trf A\trf)
\off \longrightarrow\off 
H^{\dff n}\dff(\trf X\fff,\qff A\trf)
\] 

\vspace{-8.2pt}
is defined as follows.\oss
If\qss $a\qff \in\qff H^{\dff n\dff -\dff 1}\dff (\trf A\trf)$\qss 
is\dss represented
by a cocycle $\alpha$\nnsp,\qss 
one extends $\alpha$ as a function on simplices in an arbitrary way to a cochain\qss 
$\widetilde{\alpha}\qff \in\qff C^{\dff n\dff -\dff 1}\dff (\trf X\trf)$\nnsp,\qss 
and then takes\sss the coboundary $\partial^*(\qff \widetilde{\alpha}\qff)$\nnsp.\qss 
This coboundary belongs\dss to\dss $C^{\dff n}\dff(\trf X\fff,\qff A \trf)$\nnsp.\oss
One can check\dss that\dss $\partial^*(\qff \widetilde{\alpha}\qff)$
is\dss a cocycle and\dss its\dss cohomology class 
$\delta\dff(\dff a\dff)$
does not depend on the choices of\dss $\alpha$ and\dss $\widetilde{\alpha}$\nnsp.

\myuppar{The induced\dss maps.}
A simplicial\dss map\qss
$\psi \colon\qff X\qff \longrightarrow\qff Y$\qss
leads to a map\vspace{3.8pt}
\[
\quad
\psi^* \colon\qff 
C^{\dff n}\dff(\trf Y\trf)
\off \longrightarrow\off 
C^{\dff n}\dff(\trf X\trf)
\]

\vspace{-8.2pt}
called\dss the\dss \emph{induced map}\dss and defined as follows.\qss
For an $n$\dnsp-simplex $\tau$ the cochain\dss
$\psi^*\dff(\dff \tau\trf)$\dss is defined as the formal 
sum of all $n$\dnsp-simplices of\dss $X$ mapped by $\psi$ onto $\tau$\nnsp.\oss 
The map $\psi^*$ is extended to the whole space $C^{\dff n}\dff(\trf Y\trf)$ by linearity.\qss 
The cochain-level\dss map $\psi^*$ induces a map\vspace{3.8pt}
\[
\quad
\psi^* \colon\qff 
H^{\dff n}\dff(\trf Y\dff)
\off \longrightarrow\pff 
H^{\dff n}\dff(\trf X\trf)
\]

\vspace{-8.2pt}
of cohomology groups\qss (also denoted also by $\psi^*$\nsp).\oss
The induced maps commute with\dss the connecting maps in the following sense.\oss
Let\vspace{3.8pt} 
\[
\quad
\varphi \colon (\qff X',\qff A' \trf)
\off \longrightarrow\off 
(\qff X\fff,\qff A \trf)
\]

\vspace{-8.2pt} 
be a simplicial map of pairs,\qss 
i.e.\dss a simplicial map\qss
$\varphi\dff \colon\dff X'\qff \longrightarrow\qff X$\qss 
such that\qss $\varphi\dff (\dff A'\dff )\off \subset\off A$\nnsp,\qss 
and\dss let 
$\varphi_{\dff A}\dff \colon\dff A'\qff \longrightarrow\qff A$\qss
be the map induced\dss by\sss $\varphi$\nnsp.\oss 
Then\dss the diagram\vspace{6.8pt}
\begin{equation}
\label{general-diagram}
\quad
\begin{tikzcd}[column sep=huge, row sep=mysize6]\dis
H^{\dff n\dff -\dff 1}\dff (\trf A \trf) 
\arrow[d, "\dis \varphi_{\dff A}^*"]
\arrow[r, "\dis \delta"]
& 
H^{\dff n}\dff (\qff X\fff,\qff A \trf)
\arrow[d, "\dis \varphi^*"] 
\\ 
H^{\dff n\dff -\dff 1}\dff (\trf A' \trf)
\arrow[r, "\dis \delta"]
& 
H^{\dff n}\dff (\qff X'\fff,\qff A' \trf)
\off,
\end{tikzcd}
\end{equation}

\vspace{-5.2pt}
is\dss commutative.\qss 
Clearly\halfff,\qss the diagram\qss (\ref{hdiag})\qss is a special case of\dss
the diagram\qss (\ref{general-diagram}).\oss

\myuppar{The commutativity\sss of\qss the diagram (\ref{general-diagram}).}
Let\dss us\dss outline a proof\dss of\trs this commutativity\halfff.\oss
Let\qss $a\qff \in\qff H^{\dff n\dff -\dff 1}\dff (\trf A \trf)$\dss 
be a cohomology class.\qss
It\sss can be represented\dss by\sss a cocycle\qss 
$\alpha\qff \in\qff C^{\dff n\dff -\dff 1}\dff (\trf A \trf)$\nnsp.\pss 
Let\dss us\dss extend $\alpha$\sss to a cochain\qss 
$\widetilde{\alpha}\qff \in\qff C^{\dff n\dff -\dff 1}\dff (\trf X \trf)$\qss 
as\sss in the definition of\dss the connecting map\dss $\delta$\nnsp.\qss 
Then\dss $\varphi^*\dff (\qff \partial^*(\qff \widetilde{\alpha} \qff) \qff)$\dss
is\dss a cocycle
representing\dss the cohomology class\dss
$\varphi^*\dff (\qff \delta\dff (\dff a\dff ) \qff)$\nnsp.\qss 
On\sss the other\sss hand,\qss
the cochain\dss $\varphi^*\dff (\qff \widetilde{\alpha} \qff)$\dss
extends\dss the cochain\dss $\varphi_{\dff A}^*\dff (\dff \alpha \dff)$\dss
and\dss therefore\sss the coboundary\qss
$\partial^*\dff (\qff \varphi^*\dff (\qff \widetilde{\alpha} \qff) \qff)$\dss 
represents\dss $\delta\qff (\qff \varphi_{\dff A}^*\dff(\dff a \dff) \qff)$\nnsp.\qss 
Therefore\dss
it\dss is\dss sufficient\dss to prove\sss that\vspace{2.5pt}
\begin{equation}
\label{commut-chains}
\quad
\varphi^*\dff \left(\qff \partial^*\dff (\qff \widetilde{\alpha} \qff) \qff\right)
\off =\off
\partial^*\dff \left(\qff \varphi^*\dff (\qff \widetilde{\alpha} \qff) \qff\right)\qff.
\end{equation}

\vspace{-9.5pt}
In fact\halfff,\qss
$\varphi^*\circ\pff \partial^*
\off =\off\dff 
\partial^*\circ\qff\fff \varphi^*$\nnsp,\qff\oss
i.e.\qss\vspace{2.5pt}
\begin{equation}
\label{commut-chains-general}
\quad
\varphi^*\dff \left(\qff \partial^*\dff (\qff \beta \dff) \qff\right)
\off =\off
\partial^*\dff \left(\qff \varphi^*\dff (\qff \beta \dff) \qff\right)
\end{equation}

\vspace{-9.5pt}
for every cochain\qss $\beta\qff \in\qff C^{\dff n\dff -\dff 1}\dff (\trf X \trf)$\nnsp.\qss

By\dss linearity,\oss in order\dss to prove\qss (\ref{commut-chains-general})\qss
it\dss is\dss sufficient\dss to consider\dss the case
when\dss $\beta$\dss corresponds\sss to a simplex\dss $\tau$\nnsp.\oss
In\dss this case both sides of\pss (\ref{commut-chains-general})\qss
are equal\dss to\sss the sum of\dss all $n${\dnsp}-dimensional simplices\dss 
$\sigma$\dss of\qss $X'$\dss
such\sss that\dss $\varphi\dff(\dff \sigma\dff)$\dss is\dss also
$n$\dnsp-dimensional and\dss has\dss $\tau$\dss as a\sss face.\qss 
The only subtle point\dss here\dss is\dss the fact $n$\dnsp-dimensional 
simplices $\sigma$ of\qss $X'$\dss
such that\qss $\varphi\dff(\dff \sigma\dff)\qff =\qff \tau$\qss enter into\qss
$\partial^*\dff (\qff \varphi^*\dff (\qff \beta\dff ) \qff)$\qss with the
coefficient\dss $2$\nnsp,\oss which is equal\dss to\dss $0$\dss in\dss
$\mathbf{F}_{\dff 2}$\nnsp.\oss
In\dss the context of\qss Sperner's\dss Lemma such simplices are exactly\dss the
$n${\dnsp}-dimensional simplices $\sigma$ such\dss that\dss
$\varphi\dff(\dff \sigma\dff)\off =\off \Delta_{\dff n\dff +\dff 1}$\nsp,\oss
and\dss arguments used\dss in\dss the combinatorial\dss proof\dss
work\dss in\dss the present\sss situation also.\oss

\myitpar{Cochain-level\sss proof\dss of\qss Sperner's\dss Lemma.} 
Not surprisingly,\qss 
the cochain-level proof\dss is based directly on\qss (\ref{commut-chains})\qss
or\qss (\ref{commut-chains-general}),\pss
bypassing\dss the commutativity of\dss the cohomological diagram\qss (\ref{general-diagram}).
As before,\qss we will argue by an induction by $n$\nnsp,\qss 
the case\dss $n\off =\off 0$\dss being trivial.\qss 
Suppose\sss that\dss $n\qff >\qff 0$\dss
and\dss let\sss us consider\sss the simplex 
$\Delta_{\dff n\dff +\dff 1}$ as a cochain of\qss $\partial\dff \Delta$\nnsp.\qss 
There are no $n${\dnsp}-simplices in\dss $\partial\dff \Delta$\dss
and\dss hence\dss
$C^{\dff n}\dff (\trf \partial\dff \Delta \trf)\off =\off 0$\nnsp.\qss
Therefore $\Delta_{\dff n\dff +\dff 1}$ is\dss a\sss cocycle.\oss

We would\dss like to apply\qss (\ref{commut-chains})\qss to 
$\Delta_{\dff n\dff +\dff 1}$ in the role of $\alpha$\nnsp.\qss 
In order\sss to do\sss this 
one needs to extend $\Delta_{\dff n\dff +\dff 1}$ to 
an $(\fff n\dff -\dff 1\fff)$\dnsp-dimensional cochain of\dss $\Delta$\nnsp.\oss 
Since all $(\fff n\dff -\dff 1\fff)${\dnsp}-simplices of $\Delta$ are contained in $\partial\dff \Delta$,\qss 
the only\dss possible extension of\qss $\Delta_{\dff n\dff +\dff 1}$ 
is\dss $\Delta_{\dff n\dff +\dff 1}$ itself\halfff.\oss
Hence\qss
(\ref{commut-chains})\qss takes the form\vspace{2.5pt}
\begin{equation}
\label{rc-chains}
\quad
\varphi^*\dff \left(\qff \partial^*\dff (\dff \Delta_{\dff n\dff +\dff 1}\dff ) \qff\right)
\off =\qff\off
\partial^*\dff \left(\qff \varphi^*\dff (\dff \Delta_{\dff n\dff +\dff 1}\dff ) \qff\right)
\qff.
\end{equation}

\vspace{-9.5pt}
Let us compute the left hand side of\dss the equality\qss (\ref{rc-chains}).\oss 
Clearly,\qss  
$\partial^*\dff (\trf \Delta_{\dff n\dff +\dff 1} \trf)
\off =\off \Delta$\nnsp,\qss 
and hence\vspace{2.5pt} 
\[
\quad
\varphi^*\dff \left(\qff \partial^*\dff (\qff \Delta_{\dff n\dff +\dff 1} \qff) \qff\right)
\off =\qff\off 
\varphi^*\dff \left(\qff \Delta \qff\right)\qff.
\]
 
\vspace{-9.5pt} 
By\dss the definition,\qss
$\varphi^*\dff (\trf \Delta \trf)
\off =\qff\off
\rho_{\dff 1}\qff +\qff \ldots\qff +\qff \rho_{\dff e}$\nnsp,\qss 
where\dss $\rho_{\dff i}$\dss denote\sss 
the same simplices as in\dss the combinatorial proof\halfff.\qff\oss 
It\dss follows that\vspace{2.5pt} 
\begin{equation}
\label{rhs}
\quad
\varphi^*\dff \left(\qff \partial^*\dff (\qff \Delta_{\dff n\dff +\dff 1} \qff) \qff\right)
\off =\qff\off
\rho_{\dff 1}\qff +\qff \ldots\qff +\qff \rho_{\dff e}\qff.
\end{equation}

\vspace{-9pt}
Let\sss us compute now the right hand side of\qss (\ref{rc-chains}).\oss 
The simplices\qss
$\sigma_1\dff,\pff \ldots\dff,\pff \sigma_h$\qss and\qss 
$\tau_1\dff,\pff \ldots\dff,\pff \tau_g$\qss
from\dss the combinatorial\dss proof\dss provide a complete\dss list\dss
without\dss repetitions of\dss
$(\fff n\dff -\dff 1\fff)$\dnsp-di\-men\-sion\-al simplices of\qss $S$ 
mapped\dss by\dss $\varphi$\dss onto\dss $\Delta_{\dff n\dff +\dff 1}$\nnsp.\oss
Therefore\vspace{3.25pt}
\[
\quad
\varphi^*\dff \left(\qff \Delta_{\dff n\dff +\dff 1} \qff\right)
\off =\qff\off 
\sigma_{\dff 1}\qff +\qff \ldots\dff +\dff \sigma_{\dff h}
\off +\off 
\tau_{\dff 1}\qff +\qff \ldots\qff +\qff \tau_{\dff g}\qff
\]

\vspace{-9pt}
and we need\dss to compute the coboundary\qss\vspace{3pt}
\[
\quad
\partial^*\dff (\qff
\sigma_{\dff 1}\qff +\qff \ldots\qff +\qff \sigma_{\dff h1}
\off +\off 
\tau_{\dff 1}\qff +\qff \ldots\qff +\qff \tau_{\dff g}
\qff) 
\qff.
\]

\vspace{-9pt}
Every\dss $\sigma_{\dff i}$\trs is\dss an 
$(\fff n\dff -\dff 1\fff)$\dnsp-di\-men\-sion\-al\sss simplex contained\dss in 
$\partial\dff S$\qss 
and\dss hence is a face of\dss a unique $n${\dnsp}-dimensional simplex\dss 
$\Sigma_{\dff i}$\dss of\dss $S$\nnsp.\oss 
Clearly,\qss 
$\partial^*\dff (\dff \sigma_{\dff i}\qff )
\off =\qff\off 
\Sigma_{\dff i}$.\qff\oss 
It\dss follows\sss that\vspace{4.5pt}
\begin{equation}
\label{lhs}
\quad
\partial^*\dff \left(\qff \varphi^*(\qff \Delta_{\dff n\dff +\dff 1} \qff)\qff \right)
\off =\qff\off 
\partial^*\dff \left(\qff \sigma_{\dff 1}\qff +\qff \ldots\qff +\qff \sigma_{\dff h}
\off +\off 
\tau_1\qff +\qff \ldots\qff +\qff \tau_g \qff\right)
\end{equation}

\vspace{-34.5pt}
\[
\quad
\phantom{\partial^*\dff \left(\qff \varphi^*(\qff \Delta_{\dff n\dff +\dff 1} \qff)\qff \right)
\off }
=\qff\off
\Sigma_{\dff 1}\qff +\qff \ldots\qff +\qff \Sigma_{\dff h}
\off +\off 
\partial^*\dff (\dff \tau_{\dff 1} \trf)\qff +\qff \ldots\qff +\qff \partial^*\dff (\dff \tau_{\dff g} \trf)
\qff.
\]

\vspace*{-7.5pt}
By substituting\qss (\ref{rhs})\qss and\qss (\ref{lhs})\qss into\qss (\ref{rc-chains})\qss 
we see that\vspace{3.25pt}
\begin{equation}
\label{chain}
\quad
\Sigma_{\dff 1}\qff +\qff \ldots\qff +\qff \Sigma_{\dff h}
\off +\off 
\partial^*\dff (\dff \tau_{\dff 1} \trf)\qff +\qff \ldots\qff +\qff \partial^*\dff (\dff \tau_{\dff g} \trf)
\off =\qff\off
\rho_{\dff 1}\qff +\qff \ldots\qff +\qff \rho_{\dff e}
\qff.
\end{equation}

\vspace{-8.75pt}
Every simplex\dss $\tau_k$\dss is contained in\dss the interior of\dss $S$\dss 
and hence 
is a face of\dss exactly two\dss $n${\dnsp}-dimensional simplices of\dss $S$\nnsp.\oss 
Therefore every coboundary\dss $\partial^*\dff(\dff \tau_k\dff)$\dss is a sum of\dss two simplices.\dss 
It\dss follows\sss that\dss 
the left\dss hand side of\qss (\ref{chain})\qss
is a sum of\qss
$h\qff +\qff 2\dff g$\qss simplices.\oss
Clearly,\oss the right\dss hand side of\qss (\ref{chain})\qss
is a sum of\qss $e$\dss simplices.\oss
This does not\dss implies\sss that\qss
$h\qff +\qff 2\dff g\off =\off e$\qss
because we are working over\dss $\mathbb{F}_{\dff 2}$\qss
(and\qss $h\qff +\qff 2\dff g$\qss is not equal\dss to\dss $e$\dss in general\fff).\oss
But summing\dss the coefficients of\dss simplices at\dss both sides of\qss
(\ref{chain})\qss as elements of\dss $\mathbb{F}_{\dff 2}$\dss shows that\vspace{3pt}
\begin{equation}
\label{congruence}
\quad
h\qff +\qff 2\dff g
\off \equiv\off 
e\off\off {\rm mod}\off 2
\end{equation}

\vspace{-9pt}
and\dss hence\qss
$h\off \equiv\off e\off\off {\rm mod}\off 2$\nnsp.\oss
Therefore\dss $e$\dss is odd\dss if\dss and only if\dss $h$\dss is odd.\oss  \eproof

\myuppar{The cochain-level\dss proof\dss and\dss the combinatorial\dss proof.}
The chain-level\dss proof\dss can\sss be easily\dss modified\dss
to get\sss not\sss only\dss the congruence\qss (\ref{congruence}),\pss
but\sss also\sss the equality\qss (\ref{sperner}),\oss
the heart\sss of\qss Sperner's\qss proof\halfff.\oss
As we saw,\pss the left\dss hand side of\qss (\ref{chain})\qss
is a sum of\qss
$h\qff +\qff 2\dff g$\qss simplices.\oss
Some of\trs them cancel\sss each other\halfff.\oss
An $n${\dnsp}-dimensional simplex $\sigma$ occurs in\dss this sum\dss twice\dss
if\trs either $\sigma_{\dff i}$ and\dss $\tau_{\dff k}$ are\sss faces of\dss $\sigma$\dss
for\sss some\qss $i\fff,\pff k$\qss 
({\fff}in\dss this case\qss 
$\sigma\off =\off \Sigma_{\dff i}$\nsp),\oss
or\dss $\tau_{\fff i}$\dss and\dss $\tau_{\dff k}$\dss are faces of\dss $\sigma$\dss
for some\qss $i\off \neq\off k$\nnsp.\oss
Such simplices\dss $\sigma$\dss are exactly\dss
the $n$\dnsp-dimensional\sss simplices $\sigma$ such\dss that 
$\varphi\dff(\dff \sigma\dff)\off =\off \Delta_{\dff n\dff +\dff 1}$\nsp,\oss
and\dss there are $f$ of\trs them.\oss
In other words,\oss
there are
exactly\dss $f$\dss cancellations
at\dss the left hand side of\qss (\ref{chain}).\pss
Therefore\qss (\ref{chain})\qss implies\sss that\qss
$2\dff g\dff +\dff h\qff -\qff 2\dff f\off =\off e$\qss
and\dss hence implies\qss (\ref{sperner}).\pss
One can say\dss that\qss (\ref{chain})\qss is\dss 
the cochain-level\dss realization of\trs the equation\qss (\ref{sperner}),\pss 
and\dss the cochains-based\dss proof\dss is\dss a\qss ``linearization''\qss
of\trs the combinatorial\dss proof\halfff.\oss

\mysection{Brouwer's\qss fixed-point\qss theorem}{deduction}

\myuppar{Brouwer's\qss fixed-point\dss theorem.} 
\emph{Every continuous map\qss 
$f\dff \colon\dff \Delta\qff \longrightarrow\qff \Delta$\qss has a fixed point.}

The first\dss proof\trs based on\qss Sperner's\trs lemma\dss is\dss due\sss to\qss
Knaster\halfff,\pss Kuratowski,\pss and\qss Mazurkiewich\qss \cite{kkm}\halfff.\oss 
The following proof\qss is\dss a well\dss known simplified version of\dss it\halfff.\oss
The simplification\dss results from arguing\dss by contradiction.\oss
Knaster--Kuratowski--Mazurkiewich\dss construct\sss
a labeling satisfying\dss the assumptions of\qss Sperner's\trs lemma
without\sss assuming\dss that $f$ has no fixed\dss points.

\proof 
Let us denote by $x_{\dff i}$ the $i${\dnsp}th\dss coordinate of\qss 
$x\qff \in\qff \mathbb{R}^{\dff n\dff +\dff 1}$\nnsp.\qff\oss 
We may assume that\vspace{3.875pt} 
\[
\quad
\Delta
\off =\off
\bigl\{\pff
(\dff x_{\dff 1},\pff \ldots\fff,\pff x_{\dff n\dff +\dff 1} \dff)
\pff \bigl|\pff 
x_{\dff 1},\pff \ldots\fff,\pff x_{\dff n\dff +\dff 1}
\qff \geq\qff 0
\quad\mbox{and}\quad 
x_{\dff 1}\qff +\qff \ldots\qff +\qff x_{\dff n\dff +\dff 1}
\off =\off
1
\pff\bigr\}
\off \subset\off
\mathbb{R}^{\dff n\dff +\dff 1}
\off
\]

\vspace{-8.125pt}
and\dss that\dss $\Delta_{\dff i}$\dss is\sss the face of\dss $\Delta$\dss
defined\dss by\dss the equation\qss $x_{\dff i}\off =\off 0$\nnsp.\oss 

Suppose that\qss 
$x\fff,\pff y\qff \in\qff \Delta$\qss and\dss that\qss 
$y_{\dff i}\qff \geq\qff x_{\dff i}$\qss 
for all\qss 
$i\off =\off 1\fff,\pff 2\fff,\pff \ldots\fff,\pff n\qff +\qff 1$\nnsp.\pss 
Since the coordinates of\dss $y$\dss are nonnegative and\dss 
their sum\dss is\dss equal\dss to\dss $1$\nnsp,\pss 
and\dss the same\dss is\dss true for\dss $x$\nnsp,\pss 
these inequalities imply\dss that\qss 
$y\qff =\qff x$\nnsp.\oss
It\dss follows\sss that\dss if\qss 
the map $f$ has no fixed\dss points\qss
and\dss $x\qff \in\qff \Delta$\nnsp,\pss
then\qss 
$f\dff (\dff x\dff)_{\fff i}\qff <\qff x_{\dff i}$\dss
for some\dss $i$\qss
({\fff}because\qss $f\dff (\dff x\dff)\qff \neq\qff x$\nsp).\oss
At\dss the same time,\oss if\qss
$x\qff \in\qff \Delta_{\dff i}$\nnsp,\oss
then\qss
$f\dff (\dff x\dff)_{\fff i}\pff \geq\qff x_{\dff i}$\qss
because\dss $f\dff (\dff x\dff)_{\fff i}$\dss
is always non-negative and\sss in\dss this case\qss
$x_{\dff i}\off =\off 0$\nnsp.\oss\vspace{0.5pt}

Suppose\sss that\dss $f$\dss has no fixed\dss points and
choose a sequence of subdivisions\qss 
$S_{\fff 0}\fff,\pff S_{\fff 1}\fff,\pff  S_{\fff 2}\fff,\pff \ldots$\qss 
of\dss $\Delta$\dss 
in such a way that the maximal diameter of\dss simplices of\dss 
$S_{\fff i}$\dss tends to $0$ when\qss 
$i\qff \longrightarrow\qff \infty$\nnsp.\pss 
For example,\pss one can take\qss
$S_{\fff 0}\qff =\qff \Delta$\qss and\dss 
$S_{\fff k\dff +\dff 1}$\dss 
to be the barycentric subdivision of\dss $S_{\fff k}$\nnsp.\pss\vspace{0.5pt} 

For each subdivision\dss  $S_{\fff k}$\dss and\dss each\dss vertex\dss
$w$\dss of\qss $S_{\fff k}$\dss  
let\dss us\dss label\dss $w$\dss by\dss any\dss $i$\dss such that\qss 
$f\dff(\dff w\dff)_{\fff i}\qff <\qff w_{\dff i}$\nnsp.\oss
If\qss $w\qff \in\qff \Delta_{\dff i}$\nnsp,\oss
then\dss the label of\dss $w$\dss cannot\dss be equal\dss to\dss $i$\nnsp.\oss
Hence any such\dss labeling\dss satisfies\sss the assumptions of\qss Sperner's\dss lemma.\oss 
Sperner's\dss lemma\dss implies\sss that\dss for every\dss $k$\dss 
there is a simplex\dss $\sigma_{\dff k}$\dss of\qss 
$S_{\fff k}$\dss having\qss
$\{\qff 1\fff,\pff 2\fff,\pff \ldots\fff,\pff n\qff +\qff 1 \qff\}$\qss 
as the set of\dss labels of\dss its vertices.\pss 
Let\dss $x\dff (\dff k \dff)$\dss be an arbitrary\sss point\sss of\qss $\sigma_{\dff k}$\nsp.\pss
Since\dss $\Delta$\dss is\dss compact\halfff,\pss 
one can assume\qss 
(after replacing\dss the sequence\dss $S_{\fff k}$\dss by a subsequence,\pss 
if\dss necessary)\qss 
that the sequence\dss $x\dff (\dff k \dff)$\dss 
converges\sss to some point\qss $x\qff \in\qff \Delta$\nnsp.\pss 
Since the diameters of\dss simplices\dss $\sigma_{\dff k}$\dss tend\dss to\dss $0$\nnsp,\pss 
for every\qss 
$i\off =\off 1\fff,\pff 2\fff,\pff \ldots\fff,\pff n\qff +\qff 1$\qss 
there are points\qss $w\qff \in\qff \Delta$\qss 
arbitrarily close to\dss $x$\dss and such\sss that\qss 
$f\dff(\dff w\dff)_{\fff i}\qff <\qff w_{\dff i}$\nnsp,\oss 
for example,\pss the vertices of\trs the simplices\dss $\sigma_{\dff k}$\dss labeled by\dss $i$\dss 
for sufficiently\dss big\dss $k$\nnsp.\oss
By\dss passing\dss to\dss the limit\sss we conclude that\qss 
$f\dff (\dff x\trf)_{\fff i}\pff \leq\qff x_{\dff i}$\qss
for all\qss 
$i\off =\off 1\fff,\pff 2\fff,\pff \ldots\fff,\pff n\qff +\qff 1$\nnsp.\pss 
By\dss the observation at\dss the beginning of\dss the proof\halfff,\pss 
this implies that\qss 
$f\dff(\dff x\trf)\off =\off x$\nnsp,\oss 
contrary\dss to the assumption.\oss \eproof

\myuppar{Fixed\dss point\dss free maps and\dss retractions.} 
At\dss the first sight\dss this proof\dss completely avoids a key step of\dss almost\sss 
all\dss proofs of\dss Brouwer\dss fixed-point\dss theorem\fff:\pss 
the construction of\dss a retraction\qss 
$\Delta\qff \longrightarrow\qff \partial\dff \Delta$\qss 
from a fixed\dss point\dss free map and\dss then using\dss the\qss \emph{no-retraction\dss theorem}\qss 
to the effect\dss that such\dss retractions do not exist.\pss 
Let us recall\dss this construction.\pss 

Suppose that\qss 
$f\dff \colon\dff \Delta\qff \longrightarrow\qff \Delta$\qss 
has no fixed\dss points and,\pss 
in addition,\pss that\dss the following assumption holds\fff:\oss
({\fff}A\fff)\oss 
\emph{for each\qss $x\qff \in\qff \Delta$\qss 
the segment\dss having\dss $x$\dss and\dss $f\dff(\dff x\trf)$\dss 
as its endpoints\dss is\dss not contained\dss in\dss a\qss face of\qss $\Delta$\nnsp.}\oss
Under\dss these assumptions one can define a map\vspace{4.5pt} 
\[
\quad
r\dff \colon\dff \Delta\qff \longrightarrow\qff \partial\dff \Delta
\]

\vspace{-7.5pt}
by\dss assigning to\qss $x\qff \in\qff \Delta$\qss the point\dss 
$r\dff(\dff x\trf)$\dss of\dss intersection 
of\dss the ray\dss going from\dss $f\dff (\dff x\trf)$\dss to\dss $x$\dss with\dss 
$\partial\dff \Delta$\nnsp.\oss
This ray is well-defined\dss because\qss $f\dff (\dff x\trf)\off \neq\off x$\nnsp.\pss 
The assumption\qss ({\fff}A\fff)\qss ensures that\dss the intersection 
of\qss this ray\dss with\dss $\partial\dff \Delta$\dss consists of\dss only one point
and\dss hence\dss
$r$\dss is well-defined.\pss 
Clearly,\pss $r$\sss is\dss a retraction\qss 
$\Delta\qff \longrightarrow\qff \partial\dff \Delta$\nnsp.\oss
If\qss the map\dss $f$\dss does not\sss satisfy\qss ({\fff}A\fff),\oss 
one can replace\dss $f$\dss by\sss a map\qss 
$f'\dff \colon\dff \Delta\qff \longrightarrow\qff \Delta$\qss 
which\dss does 
and still\dss has no fixed\dss points.\pss 
In fact,\pss 
if\qss 
$f'\dff(\trf \Delta\trf)\qff \subset\qff \Delta\qff \smallsetminus\qff \partial\dff \Delta$\nnsp,\pss
then\dss $f'$\dss satisfies\qss ({\fff}A\fff).\pss 
Clearly\halfff,\pss such a map\dss $f'$\dss can be chosen to be arbitrarily close to\dss $f$\dnsp,\pss
and if\dss $f'$\dss is\dss sufficiently close to\dss $f$\dnsp,\oss 
it\dss has no fixed points together\dss with\dss $f$\nnsp.\pss 

The standard expositions usually deal with\dss maps from a ball\dss to itself\halfff.\pss 
Since a ball\dss is\dss strictly convex in
contrast\dss with\dss $\Delta$\nnsp,\pss 
the analogue of\qss ({\fff}A\fff)\qss is\dss automatically satisfied.\pss 
This suggests
another construction of\dss $r$\nnsp.\pss 
Let\qss 
$h\dff \colon\dff \Delta\qff \longrightarrow\qff B$\qss 
be a homeomorphism.\pss
Then\qss 
$h\dff\circ f \circ\dff h^{\dff {}-\dff 1}$\qss
is\dss a continuous map\qss  
$B\qff \longrightarrow\qff B$\qss
without\dss fixed\dss points
and\dss hence\sss leads\sss to a\sss retraction\qss
$r_{\qff B}\dff \colon\dff B\qff \longrightarrow\qff \partial\dff B$\nnsp.\qss
The map\qss
$r
\off =\off
h^{\dff {}-\dff 1}\circ\dff r_{\qff B}\dff \circ\dff h$\qss
is\dss a\sss retraction\qss
$\Delta\qff \longrightarrow\qff \partial\dff \Delta$\nnsp.\oss

\myuppar{The hidden retraction.}
In fact,\pss the retraction\sss $r$\sss is\dss implicitly\dss used in\dss 
Knaster-Kuratowski-Mazurkiewich\dss proof\halfff.\pss 
To simplify\dss the discussion,\pss let\sss us assume that\dss $f$\dss satisfies\sss
the assumption\qss ({\fff}A\fff).\oss
Let\dss $w$\dss be a vertex.\oss
The vector\dss
$w\qff -\qff f\dff(\dff w\trf)$\dss
is\dss parallel\dss to the hyperplane\qss\vspace{4.5pt} 
\[
\quad
x_{\dff 1}\qff +\qff \ldots\qff +\qff x_{\dff n\dff +\dff 1}
\off =\off
1
\qff
\] 

\vspace{-7.5pt}
and\dss hence\sss the ray\sss from $f\dff(\dff w\trf)$ 
to $w$ is\sss contained\sss in\sss this hyperplane 
and intersects\dss $\partial\dff \Delta$\nnsp.\pss 
The vertex $w$ may\dss be labeled\dss by\sss $i$\trs 
if\trs and\dss only\dss if\qss
$w_{\dff i}\qff >\qff f\dff(\dff w\dff)_{\fff i}$\nnsp,\oss
i.e.\qss if\trs and\dss only\dss if\qss
$w\qff -\qff f\dff(\dff w\trf)$\qss 
has positive\dss $i${\dnsp}th coordinate.\oss 
By\dss the construction of\dss the retraction\dss $r$\dss
this condition is equivalent\dss to\qss
$r\dff(\dff w\trf)\qff -\qff w$\qss having\sss positive $i${\dnsp}th\dss coordinate.\oss
Hence $w$ may\dss be labeled\dss by\dss $i$\qss 
if\trs and\dss only\trs if\dss 
$r\dff(\dff w\trf)$\dss 
is contained in\dss the intersection of\qss 
$\partial\dff \Delta$\dss 
with\dss the open half-space\qss\vspace{4.5pt} 
\[
\quad
\bigl\{\pff
(\qff x_{\dff 1},\pff \ldots\fff,\pff x_{\dff n\dff +\dff 1} \qff)
\pff \bigl|\pff 
x_{\dff i}\qff >\qff w_{\fff i}
\pff\bigr\}
\off.
\] 

\vspace{-7.5pt}
This\dss intersection\dss is\dss contained\dss in\qss 
$\partial\dff \Delta\qff \smallsetminus\qff \Delta_{\dff i}$\qss 
because\sss the face\dss $\Delta_{\dff i}$\dss is\dss defined\dss by\qss 
$x_{\dff i}\off =\off 0$\nnsp.\oss
It\dss is\dss tempting\dss to allow as a\dss label of\dss $w$\dss any\dss $i$\dss 
such that\qss 
$r\dff(\dff w\trf)\qff \in\pff \partial\dff \Delta\qff \smallsetminus\qff \Delta_{\dff i}$\nnsp.\oss 
This condition\dss has the advantage of\qss 
$\Delta\qff \smallsetminus\qff \Delta_{\dff i}$\qss 
being\dss independent of\dss $w$\nnsp,\oss 
but\sss after passing\dss to\sss the limit\dss the sets\qss 
$\Delta\qff \smallsetminus\qff \Delta_{\dff i}$\qss
has\sss to be replaced\dss by\dss their closures,\oss
which,\oss
in contrast\sss with\dss the sets\qss 
$\Delta\qff \smallsetminus\qff \Delta_{\dff i}$\qss
themselves,\oss
have non-empty\sss intersection\qss
(for example,\oss every vertex of\dss $\Delta$\dss belongs\sss to\sss their\sss intersection).\oss
But\sss one can avoid\dss passing\dss to\sss the limit\dss
by\dss using compactness in a less direct\dss manner
and strengthening\dss the condition\qss
$r\dff(\dff w\trf)\qff \in\qff \partial\dff \Delta\qff \smallsetminus\qff \Delta_{\dff i}$\qss
in a standard\dss in combinatorial\dss topology\sss way\halfff.\oss
This leads\sss to a proof\sss of\dss 
the no-retraction\dss theorem\dss based\sss on\qss Sperner's\trs lemma.\oss

\myuppar{No-retraction\dss theorem.} 
\emph{There exists no retraction\qss $\Delta\qff \longrightarrow\qff \partial\dff \Delta$\nnsp.}

\proof 
Suppose that\dss $r$\dss is a such retraction.\oss 
Recall\dss that\dss the\qss \emph{open star}\qss 
${\rm st}\trf(\dff w\trf)$\dss of\dss a vertex\dss $w$\dss 
of\dss a simplicial complex is the union 
of\dss all simplices having\dss $w$\dss as a vertex\dss 
with\dss the faces opposite of\dss $w$\dss removed.\oss
In\dss particular,\qss 
$\partial\dff \Delta\dff \smallsetminus\dff \Delta_{\dff i}
\off =\off
{\rm st}\trf(\dff v_i\trf)$\nnsp,\pss
where the star is taken in\dss $\partial\dff \Delta$\nnsp.\oss 
By\dss the well\dss known\dss Lebesgue\dss lemma applied\dss to the open covering of\dss 
$\Delta$\dss by\dss the preimages\qss 
$r^{\dff {}-\dff 1}\dff (\qff {\rm st}\trf (\dff v_i\trf)\qff)$
every subset\sss of\sss $\Delta$ having sufficiently small diameter is
contained in one of\dss these preimages.\oss 
Applying\trs Lebesgue\trs lemma\dss is\dss simply\sss another\sss way\dss to use\sss 
the compactness of\dss $\Delta$\nnsp.\oss
Clearly,\pss
the diameter of\dss a star of\dss a vertex of\dss a simplicial\sss complex\sss 
is no bigger\dss than\dss twice 
the maximum\sss of\dss the diameters of\dss simplices of\dss this complex.\oss 
Therefore,\pss if\qss diameters of\dss simplices of\dss 
a subdivision\dss $S$\dss of\dss $\Delta$\dss are small\sss enough,\pss 
then for every vertex $w$ of\dss $S$\sss the image 
$r\dff(\qff {\rm st}\trf(\dff w\trf)\qff)$ is\dss contained in some star 
${\rm st}\trf(\dff v_i\trf)$\nnsp.\oss 
Let\dss us label each vertex\dss $w$\dss by any\dss $i$\dss such that\vspace{3pt}  
\[
\quad
r\dff(\qff {\rm st}\dff (\dff w\trf)\qff)
\off \subset\qff\off
\partial\dff\Delta
\qff \smallsetminus\qff 
\Delta_{\dff i}
\off =\off
{\rm st}\trf (\dff v_i \trf)
\qff.
\] 

\vspace{-9pt}
By\dss Sperner's\dss lemma
there\dss is\dss a\sss simplex $\sigma$ of\dss $S$ having\dss
$\{\qff 1\fff,\pff 2\fff,\pff \ldots\fff,\pff n\qff +\qff 1 \qff\}$\dss 
as\sss the set\sss of\dss labels of\dss its\sss vertices.\oss 
The interior\dss ${\rm int}\qff \sigma$\dss of\dss $\sigma$\dss 
is\dss contained\dss in\dss 
${\rm st}\trf (\dff w\trf)$\dss 
for every\sss vertex $w$ of\dss $\sigma$\nnsp.\oss 
Since for each\qss 
$i\off =\off 1\fff,\pff 2\fff,\pff \ldots\fff,\pff n\qff +\qff 1$\qss 
some vertex $w$ of\dss $\sigma$\dss is labeled\dss by\dss $i$\nnsp,\oss 
it\dss follows\dss that\qss \vspace{3pt}
\[
\quad
r\dff (\qff {\rm int}\qff \sigma \qff)
\off \subset\off
{\rm st}\trf(\dff v_i \trf)
\] 

\vspace{-9pt}
for every $i$\nnsp.\oss 
But the sets\qss 
${\rm st}\trf(\dff v_i \trf)
\off =\off
\partial\dff \Delta\qff \smallsetminus\qff \Delta_{\dff i}$\qss 
obviously\dss have empty\dss intersection.\oss 
The contradiction between\dss the\sss last\dss two statements completes\sss the proof.\oss \eproof

\myuppar{Simplicial approximations.} 
The\sss labeling\dss constructed\dss in\dss the above\sss proof\qss 
leads a\qss simplicial\dss map\qss %\emph{simplicial approximation}\qss
$\varphi\dff \colon\dff S\qff \longrightarrow\qff \partial\dff \Delta$\qss 
which\dss is\dss a\qss \emph{simplicial approximation}\qss 
of\dss the contunuous map $r$\nnsp.\oss 
In\dss fact\halfff,\pss the best\dss way\dss to define\qss 
\emph{simplicial\sss approximations}\qss is\dss to require\sss that\dss
the condition\vspace{3pt} 
\[
\quad
r\dff(\qff {\rm st}\trf (\dff w\trf) \qff)
\off \subset\off 
{\rm st}\trf (\qff \varphi\dff (\dff w\trf) \qff)
\] 

\vspace{-9pt}
holds for every\dss vertex $w$\nnsp.\oss
Nowadays usually\sss another definition\dss is\dss adopted,\oss
followed\dss by\sss a proof\trs that\dss this condition\dss is\dss
necessary and sufficient\dss 
for\dss $\varphi$\dss to be a simplicial approximation 
of\dss the continuous map\sss $r$\nnsp.\oss 
The strengthening of\trs the condition\qss
$r\dff(\dff w\trf)\qff \in\pff \partial\dff \Delta\qff \smallsetminus\qff \Delta_{\dff i}$\qss
to\vspace{3pt}
\[
\quad
r\dff(\qff {\rm st}\dff (\dff w\trf) \qff)
\off \subset\qff\off
\partial\dff\Delta
\qff \smallsetminus\qff 
\Delta_{\dff i}
\]

\vspace{-9pt}
was motivated\sss exactly\dss by\dss this\sss property\sss of\dss simplicial\sss approximations.\oss 
The use of\qss Lebesgue\dss lemma above\dss is\dss nothing else but\dss 
the standard\dss way\dss to establish\dss the existence of\dss simplicial approximations.\oss 
Sperner's\qss lemma\sss implies\sss that\dss the existence of\dss simplicial\sss approximations
of\dss a\sss retraction\qss $\Delta\qff \longrightarrow\qff \partial\dff \Delta$\qss
leads\sss to a contradiction.\oss

\vspace{6ex}
\begin{flushright}

June\qss 29,\oss 2009\oss (the first\sss version)

August\pss 24,\oss 2019\oss (the current\sss version)
 
https:/\!/\hspace*{-0.07em}nikolaivivanov.com

E-mail\halfff:\oss nikolai.v.ivanov{\fff}@{\dff}icloud.com
\end{flushright}

\end{document}